\documentstyle{article}

         \def\ba{\begin{array}}
         \def\ea{\end{array}}
         \def\be{\begin{equation}}
         \def\ee{\end{equation}}
         \def\sl{\hbox{{\rm sl}}}
         \def\half{{1 \over 2}}
         \def\Uq#1{{\rm U}_q \left( #1 \right) }
         \def\Uh#1{{\rm U}_h \left( #1 \right) }
         \def\Um#1{{\rm U}_{\mu} \left( #1 \right) }
         \def\U#1{{\rm U} \left( #1 \right) }
         \def\R{{\rm R}}
         \def\Real{\hbox{{\rm R}}}
         \def\so{{\rm so}}
         \def\poi#1{\hbox{{\rm p}}(#1)}
         \def\dd#1{{\partial\over{\partial#1}}}
\begin{document}
\hfill{}

\vbox{
    \halign{#\hfil         \cr
            q-alg/9512022 \cr\noalign{\vskip -0.5cm}
           } 
      }  
\vskip 30 mm
\leftline{ \Large \bf
          The universal R matrix for the Jordanian deformation of}
\leftline{ \Large \bf
          sl(2), and the contracted forms of so(4)}
\vskip 10 mm
\rm
\leftline{ \bf
          A. Shariati ${}^{1,2,*}$,
          A. Aghamohammadi ${}^{1,3}$
          M. Khorrami ${}^{1,2,4}$,
          }
\vskip 10 mm
{\it
  \leftline{ $^1$ Institute for Studies in Theoretical Physics and
            Mathematics, P.O.Box  5746, Tehran 19395, Iran. }
  \leftline{ $^2$ Institute for Advanced Studies in Basic Sciences,
             P.O.Box 159, Gava Zang, Zanjan 45195, Iran. }
  \leftline{ $^3$ Department of Physics, Alzahra University,
             Tehran 19834, Iran. }
  \leftline{ $^4$ Department of Physics, Tehran University,
             North-Kargar Ave. Tehran, Iran. }

  \leftline{ $^*$ E-Mail: shahram@irearn.bitnet}
  }
\begin{abstract}
We introduce a universal R matrix for the Jordanian deformation of
$\U{ \sl(2)}$. Using
$\Uh{\so(4)}=\Uh{\sl(2)} \oplus {\rm U}_{-h}(\sl(2))$, we obtain
the universal R matrix for $\Uh{\so(4)}$.
Applying the graded contractions on the universal R matrix of
$\Uh{\so(4)}$, we show that there exist three distinct
R matrices for all of the contracted algebras. It is shown that
$\Uh{\sl(2)}$, $\Uh{\so(4)}$, and all of these contracted algebras are
triangular.

\end{abstract}
\vskip 10 mm
\section{ Introduction }
The group SL(2) admits two distinct quantizations, one is the well known
Drinfeld-Jimbo ($q$-) deformation, and the other is the so called Jordanian
($h$-) deformation \cite{Manin,Ohn}. In fact the $h$-deformation itself
can be obtained by a contraction procedure from the $q$- deformation
\cite{AKS}. So far, however, there has been no expression for the Universal
R matrix of $\Uh{\sl(2)}$. It is worthy of mention that the R matrix which
was introduced in \cite{Ohn} does not satisfy the quantum Yang--Baxter
equation \cite{Vlad}. The universal R matrix for the positive Borel
subalgebra of the Jordanian $\sl(2)$ was introduced in \cite{Vlad}.

Another interesting problem is the study of non semisimple quantum
groups. One of the techniques for constructing inhomogeneous quantum
groups is contraction. There are two distinct deformations of $\U{\rm e(2)}$
both of which can be obtained by contraction of $\Uq{\rm su(2)}$,
neither of them has a universal R matrix \cite{Vaks,CGST,BCGST}.
There is also a deformation of the two dimensional Poincar\'e
algebra $\Um{\poi{2}}$ which was obtained by a contraction of
$\Uh{\sl(2,\Real)}$ and has a universal R matrix \cite{KSAA}.
In \cite{BHOS} two copies of the Jordanian deformation of $\sl(2)$ have been
used to construct the deformed algebra of $\so(4)$. Then, the process of
graded contraction \cite{DM} has been used to construct a deformation for a
fairly large class of non semisimple algebras.

In this article, we first introduce an expression for the universal R
matrix of $\Uh{\sl(2)}$, and show that this algebra is triangular
\cite{Majid}. In fact we prove that the universal R matrix
obtained in \cite{Vlad} for the positove Borel subalgebra of the Jordanian
sl(2), is the universal R matrix for the whole of $\Uh{\sl(2)}$.
Then we complete the study of \cite{BHOS}, that is, we list all possible
contractions of $\Uh{\so(4)}$ and, using the R matrix of $\Uh{\sl(2)}$,
we obtain the R matrices of $\Uh{\so(4)}$, and all its contracted algebras.
As we will see, there are three distinct R matrices for all of the contracted
algebras. It is also seen that $\Uh{\so(4)}$, and all of these contracted
algebras are triangular.

\section{ The Universal R matrix}
A quasitriangular Hopf algebra is a Hopf algebra with a universal R matrix
satisfying
\be \label{quasitri1}
\ba{ll}
( \Delta \otimes 1 ) R = R_{13} R_{23}, &
( 1 \otimes \Delta ) R = R_{13} R_{12},  \cr
\ea \ee
and
\be \label{quasitri2}
R \Delta (\cdot ) R^{-1} = \Delta' ( \cdot ) := \sigma \circ \Delta (\cdot ).
\ee
where $\sigma$ is the flip map: $\sigma (a \otimes b) = b \otimes a $.
If in addition
\be \label{triangular}
\sigma ( R^{-1} ) = R,
\ee
the Hopf algebra is called  triangular \cite{Majid}.

The Jordanian deformation of sl(2) is defined \cite{Ohn} through
\be\ba{l}
[J^3,J^+]=2{{\sinh (hJ^+)}\over h}\cr [J^3,J^-]=-\big[\cosh (hJ^+)J^- +
J^-\cosh (hJ^+)\big]\cr [J^+,J^-]=J^3,\cr\ea\ee
and
\be\ba{ll}
\Delta J^+=J^+\otimes 1+1\otimes J^+ & \Delta J^i=J^i\otimes e^{hJ^+}+
e^{-hJ^+}\otimes J^i\cr \epsilon (X)=0 & \gamma (X)=-e^{hJ^+}Xe^{-hJ^+}\cr
{\rm where}\; i=-,3,\;{\rm and}\; X\in\{ J^+,J^-,J^3\} .\cr\ea\ee
Vladimirov \cite{Vlad} has considered the subalgebra of $\Uh{\sl(2)}$
generated by the two generators $J^+$ and $J^3$ and has found the following
universal R matrix for this subalgebra,
\be \label{RVlad}
\R = \exp \left\{ { { \Delta (hJ^+)} \over \sinh ( \Delta (hJ^+) ) }
\left[ J^3 \otimes \sinh ( h J^+) - \sinh (h J^+)
\otimes J^3 \right] \right\}.
\ee
This means that R satisfies (\ref{quasitri1}) and also,
\be
R\ \Delta J^+\ R^{-1} = \Delta' J^+ ,\qquad
R\ \Delta J^3\ R^{-1} = \Delta' J^3.
\ee
We are going to show that the above  expression for R is,
in fact, the universal R matrix for the whole algebra
$\Uh{\sl(2)}$. To do so, we must show that
\be \label{jminus}
\R \Delta J^- \R^{-1} = \Delta' J^-.
\ee
Now we define E through
\be \label{E}
\R \Delta J^- \R^{-1} := \Delta' J^- + E.
\ee
To prove (\ref{jminus}), is equivalent to prove $E=0$.
{}From the commutation relations, it is seen that the power of $J^-$ in
the right hand side of (\ref{E}) does not exceed one. So,
using Hausdorf identity it can be shown that
\be
E = ( J^- \otimes 1 ) C(J^+_1, J^+_2, J^3_1, J^3_2)
+ (1 \otimes J^- ) D(J^+_1, J^+_2, J^3_1, J^3_2)
+ F(J^+_1, J^+_2, J^3_1, J^3_2).
\ee
where $ A_1 := A \otimes 1 $, and $ A_2 := 1 \otimes A $.
The first step is to show that $ C = D = 0 $. To do so, we use the fact that
the matrix (\ref{RVlad}) is a universal R matrix for the contracted form of
$\Uh{\sl(2)}$ \cite{KSAA}. The contraction procedure can be achieved
through the definitions
\be \ba{lll} P^- := \epsilon J^-, & P^+ := J^+, & J^3 := J^3. \cr \ea \ee
Note that, in this basis it is not neccessary to redefine both $J^-$ and
$J^+$. This appears slightly different from the procedure which we
introduced in \cite{KSAA}, but it is easy to see that the resulting Hopf
algebra is the same.

Now we use the above redefinitions in (\ref{E}). Note that, as $\R$ is a
function of only $J^+$ and $J^3$, it does not change in this procedure.
The result is that
\be
\R \Delta P^- \R^{-1} := \Delta' P^- + \lim_{\epsilon \to 0}\epsilon E.
\ee
But in \cite{KSAA} we have shown that
\be
\R \Delta P^- \R^{-1} := \Delta' P^-.
\ee
So, it is deduced that $\lim_{\epsilon \to 0} \epsilon E = 0 $, in which
we have to express the relation in terms of the new generators $P^{\pm}$
and $J^3$. This results in
\be (P^- \otimes 1) C + (1 \otimes P^- ) D = 0, \ee
which means that $ C = D = 0 $. We have ruled out the $J^-$ dependence
in $E$. Now we can rewrite $E$ as
\be \label{iii}
E = \sum_{m,n \geq 0} ( J^3 \otimes 1 )^m B^n g_{m,n}(J^+_1,J^+_2), \ee
where
\be B:= \half [J^3 \otimes \sinh (h J^+) - \sinh (h J^+) \otimes J^3 ]. \ee
Note that the commutation relations permit us to write $E$ as (\ref{iii}).
Using the commutation relations of $J^+$ and $J^3$ with $J^-$, the fact that
$\Delta$ and $\Delta'$ are homomorphisms of the algebra, and that
$\Delta' J^+ = \Delta J^+$, one can see that
\be \label{iv} [ \Delta J^+ , E ] = 0, \ee
and
\be \label{v} [\Delta' J^3 , E ] =
             - \{ E \cosh ( h \Delta J^+) + \cosh ( h \Delta J^+) E \}.
\ee
{}From (\ref{iv}) and  $ [ \Delta J^+ , B ] = 0 $, it is seen that
$ g_{m,n} = 0$ if $m > 0$. So,
\be \label{vi}
E = \sum_{n \geq 0} B^n g_n(J^+_1, J^+_2). \ee
Now, $E$ is an analytic function of $h$. Suppose that the smallest power of
$h$ in $E$ is $m$: $E = h^m E_{m} + O(h^{m+1})$.
Then we can deduce from (\ref{v}) that
\be \lim_{h \to 0} [ \Delta' J^3 , h^{-m} E ] = - \lim_{h \to 0}
    \{ h^{-m} E \cosh ( h \Delta J^+) + \cosh ( h \Delta J^+) h^{-m} E \},
\ee
or,
\be \label{vii}
\lim_{h \to 0} [ 1 \otimes J^3 + J^3 \otimes 1 , E_m ] = - 2 E_m. \ee
{}From (\ref{vi}), it is seen that
\be
E_m = \sum_{n \geq 0} \tilde{B}^n \tilde{g}_n (J^+_1, J^+_2),
\ee
where $\tilde{B}:=\half J^3 \otimes J^+ - J^+ \otimes J^3 $, and
$\tilde{g}_n := \lim_{h \to 0} h^{n -m} g_n$. It is easy to see that
\be \lim_{h \to 0} [ 1 \otimes J^3 + J^3 \otimes 1, \tilde{g}_n
( J^+_1 ,J^+_2 )] = 2 ( J^+_1 \dd{J^+_1} + J^+_2 \dd{J^+_2} ) \tilde{g}_n
( J^+_1, J^+_2). \ee
One can then write (\ref{vii}) as
\be \sum_{n} \tilde{B}^n ( n+1 +J^+_1 \dd{J^+_1} +
J^+_2 \dd{J^+_2} ){\tilde g}_n=0,\ee
which gives
\be  ( n+1 +J^+_1 \dd{J^+_1} + J^+_2 \dd{J^+_2} )
\tilde{g}_n = 0. \ee
This means that $ \tilde{g}_n$ is a homogeneous function of order $-(n+1)$
of $J^+_1$ and $J^+_2$. From this it is concluded that
$ \tilde{B}^n \tilde{g}_n$ is a homogeneous function of order $-1$ of
$J^+_1$ and $J^+_2$, which means that $E_m$ is such a function. But this
is impossible for $E_m \neq 0$, because $E$, and hence $E_m$, is an
analytic fuction of $J^+_i$ and $J^3_i$ (This can be deduced from the
analyticity of R and the commutation relations), and there exist no
analytic function which is homogeneous of a negative order. So, $E_m$
should be zero. This means that $E$ is zero, because we assumed that the
lowest order term of $E$ is $h^m E_m$.

So (\ref{jminus}) is correct. This completes the proof. It is worthy
of mention that (\ref{RVlad}) can also be written  in the form
\be \label{RVlad2}
\R = \exp \left\{ {\Delta -\Delta'\over 2}\left[
J^3 {hJ^+\over \sinh hJ^+}\right] \right\}.
\ee
As R is of the form $\R =\exp [(\Delta -\Delta')X]$, it is obvious that the
algebra $\Uh{\sl(2)}$ is triangular. In fact, as it will be seen, the
R matrix of $\Uh{\so(4)}$ and all its contractions which we consider, have
this property. So, all of these algebras are triangular.

In ref. \cite{BHOS}, starting from the Jordanian deformation
$\Uh{\sl(2)}$, and using
\be
\Uh{\so(4)}= \Uh{\sl(2)} \oplus {\rm U}_{-h}(\sl(2)),
\ee
$\Uh{\so(4)} $ has been constructed. Consider $\Uh{\sl(2)}$
with the generators $\{ J^3_1, J^{\pm}_1  \}$, and ${\rm U}_{-h}(\sl(2))$
with the generators $\{ J^3_2, J^{\pm}_2  \}$.The set of generators
$\{J^3,J^{\pm},N^3,N^{\pm}\}$ defined by
\be
  J^i:=J^i_1+J^i_2 \qquad N^i:=N^i_1+N^i_2  ,\qquad i=+,-,3
\ee
closes the Hopf algebra $\Uh{so(4)}$, with the following Hopf structure.
\be \label{x} \ba{l}
[J^3,J^+]={4\over h} \sinh({h\over 2} J^+)\cosh({h\over 2} N^+)\cr
[J^3,J^-]=-J^-\cosh({h\over 2} J^+)\cosh({h\over 2} N^+)-
         \cosh({h\over 2} J^+)\cosh({h\over 2} N^+)J^-\cr
\hskip 1.5cm           -N^-\sinh({h\over 2} J^+)\sinh({h\over 2} N^+)-
          \sinh({h\over 2} J^+)\sinh({h\over 2} N^+)N^-\cr
[J^3,N^+]={4\over h} \sinh({h\over 2} N^+)\cosh({h\over 2} J^+)\cr
[J^3,N^-]=-N^-\cosh({h\over 2} J^+)\cosh({h\over 2} N^+)-
          \cosh({h\over 2} J^+)\cosh({h\over 2} N^+)N^-\cr
\hskip 1.5cm           -J^-\sinh({h\over 2} J^+)\sinh({h\over 2} N^+)-
          \sinh({h\over 2} J^+)\sinh({h\over 2} N^+)J^-\cr
[N^3,N^{\pm}]=[J^3,J^{\pm}],\qquad [N^3,J^{\pm}]=[J^3,N^{\pm}]\cr
[J^+,J^-]=[N^+,N^-]=J^3,\qquad [J^{\pm},N^{\mp}]=\pm N^3,\qquad [J^i,N^i]= 0
                        \qquad {\rm where\;\;} i=+,-,3\cr
     \ea
\ee
and,
\be \label{xx}  \ba{l}
\Delta J^+  =  1 \otimes J^+ + J^+ \otimes 1 \cr
\Delta N^+  =  1 \otimes N^+ + N^+ \otimes 1 \cr
\Delta J^i  =  e^{-\frac{h}{2}N^+} \cosh(\frac{h}{2} J^+) \otimes J^i +
J^i \otimes \cosh(\frac{h}{2} J^+)e^{\frac{h}{2}N^+}  \cr
\hskip 1cm -e^{-\frac{h}{2}N^+}\sinh(\frac{h}{2} J^+) \otimes N^i
+ N^i \otimes \sinh(\frac{h}{2} J^+) e^{\frac{h}{2}N^+} \cr
\Delta N^i  =  e^{-\frac{h}{2}N^+} \cosh(\frac{h}{2} J^+)\otimes N^i
+N^i \otimes \cosh(\frac{h}{2} J^+)e^{\frac{h}{2}N^+} \cr
\hskip 1cm -e^{-\frac{h}{2}N^+} \sinh(\frac{h}{2} J^+) \otimes J^i
+J^i \otimes \sinh(\frac{h}{2} J^+)e^{\frac{h}{2}N^+}\qquad{\rm where\;\;}
                               i=-,3 \cr
\epsilon (X) = 0, \quad \gamma(X)=-e^{hN^+}Xe^{-hN^+}, \quad
{\rm where\;\;} X \in \{ J^3, J^{\pm}, N^3, N^{\pm} \} .\cr \ea \ee
This algebra has a universal R matrix which is simply the product
of the two copies of (\ref{RVlad}). The resulting R matrix, in terms of
the new generators, is
\be \label{Rso}
\R = \exp \left\{  \frac{h}{2}
(\Delta- \Delta'){\left [(J^3J^++N^3N^+)
\sinh {hJ^+\over 2}\cosh {hN^+\over 2}
-( J^3N^++N^3J^+) \sinh { h N^+\over 2}\cosh {hJ^+\over 2}\right ]\over
\cosh hJ^+-\cosh hN^+}  \right\} .
\ee

The definitions of the classical gradations and contractions have been
extended to the quantum case, by assuming that they act on the algebra of the
generators as in the classical case, and that $e^{-\frac{h}{2}N^+}$
is invariant under the quantum gradations and contractions \cite{BHOS}.
\be
(J^3, J^{\pm}, N^{\pm}, N^3, h) = (\hat J^3,
\frac{\hat J^{\pm}}{\sqrt{\mu_2 \mu_3}},
\frac{\hat N^{\pm}}{\sqrt{\mu_1 \mu_2}},
\frac{\hat N^3}{\sqrt{\mu_1 \mu_3}},
 \sqrt{\mu_1 \mu_2} \hat h ),
\ee
where $\mu_i$ fall into the one of the values $\pm 1, 0$. (We use
complex algebras, so $\mu_i$ is 1, or 0.) As it is shown
in \cite{BHOS}, for $\mu_3 \to 0$, the resulting Hopf algebra is not
well defined. However, one can modify the above contraction such that it
is well behaved for $\mu_3 \to 0$ too.
The contraction is as follows.
\be
(J^3, J^{\pm}, N^{\pm}, N^3, h) = (\hat J^3,
\frac{\hat J^{\pm}}{\sqrt{\mu_2 \mu_3}},
\frac{\hat N^{\pm}}{\sqrt{\mu_1 \mu_2}},
\frac{\hat N^3}{\sqrt{\mu_1 \mu_3}},
 \mu_3 \sqrt{\mu_1 \mu_2} \hat h ),
\ee
Applying this contraction to (\ref{x},\ref{xx}), everything
remains well behaved.
The same is true for the universal R matrix. In this way one obtaines
three distinct R matrices which we present here. The full Hopf structure
of the contracted algebras are given in the appendix.
\begin{enumerate}
\item For $\mu_3 = 0$ ; the algebras become classic (non deformed),
although the coproducts remain deformed, and
\be \R = \exp \left( \frac{\Delta - \Delta'}{2}\hat J^3 \right). \ee
\item For $\mu_3=1$, $\mu_1=0$,
\be
\R = \exp \left\{  \frac{\hat h}{2}
(\Delta- \Delta') \{{ [{\hat h\over 2}\hat N^3 \hat N^+ \hat J^+
\cosh {\hat h \hat N^+\over 2}
-( \hat J^3 \hat N^+ + \hat N^3 \hat J^+) \sinh ({ \hat h \hat N^+
\over 2})] \over 1-\cosh \hat h \hat N^+ }  \} \right\}.
\ee
\item For $\mu_3=\mu_1=1$, $\mu_2=0$,
\be
\R = \exp \left\{  \frac{\hat h}{2}
(\Delta- \Delta'){\left [(\hat J^3 \hat J^+ + \hat N^3 \hat N^+)
\sinh {\hat h \hat J^+ \over 2} \cosh {\hat h \hat N^+\over 2}
-( \hat J^3 \hat N^+ + \hat N^3 \hat J^+) \sinh { \hat h \hat N^+\over 2}
\cosh {\hat h \hat J^+ \over 2}\right ]\over
\cosh \hat h \hat J^+ - \cosh \hat h \hat N^+}  \right\}.
\ee
\end{enumerate}
In fact the last R matrix is the R matrix for
${\rm U}_{\hat h}({\rm iso(2)})\oplus {\rm U}_{-\hat h}(\rm iso(2))$,
and it is simply the product of two copies of the universal
R matrix of the deformed $\U{\rm iso(2)}$.
Note that, by $\mu_i=0$, it is meant that one must set $\mu_i=\epsilon$, and
obtain the Hopf structure in the limit $\epsilon\to 0$.

\section{Appendix}
Here we present the full Hopf structure of the contracted algebras; that is
the commutation relationships, and the nontrivial coproducts and antipodes.
Throughout this appendix, $i=+,-,3$, $j=-,3$, and $X \in \{ \hat J^3,
\hat J^{\pm}, \hat N^3, \hat N^{\pm} \}$.
\begin{enumerate}
\item $ ( \mu_1 , \mu_2 , \mu_3 ) = ( 1 , 1 , 0 ) $
$$\ba{ll}
[\hat J^3,\hat J^\pm ]=\pm 2\hat J^\pm & [\hat J^3,\hat N^\pm ]=\pm 2
\hat N^\pm\cr [\hat N^3,\hat N^\pm ]=\pm 2\hat J^\pm & [\hat N^3,\hat J^\pm ]
=0\cr [\hat J^+,\hat J^-]=0&[\hat N^+,\hat N^-]=\hat J^3\cr
[\hat J^\pm,\hat N^\mp]=\pm \hat N^3 & [\hat J^i,\hat N^i]=0\cr
  \ea$$
$$\Delta \hat J^3=1\otimes \hat J^3+\hat J^3\otimes 1+{\hat h\over 2}
(\hat N^3\otimes \hat J^+-\hat J^+\otimes \hat N^3)$$
$$\Delta \hat N^-=1\otimes \hat N^-+\hat N^-\otimes 1+{\hat h\over 2}
(\hat J^-\otimes \hat J^+-\hat J^+\otimes \hat J^-)$$
This is a deformation of U(iso(3)).
\item $ ( \mu_1 , \mu_2 , \mu_3 ) = ( 1 , 0 , 0 ) $
$$\ba{ll}
[\hat J^3,\hat J^\pm ]=\pm 2\hat J^\pm & [\hat J^3,\hat N^\pm ]=\pm 2
\hat N^\pm\cr [\hat N^3,\hat N^\pm ]=\pm 2\hat J^\pm & [\hat N^3,\hat J^\pm ]
=0\cr [\hat J^+,\hat J^-]=[\hat N^+,\hat N^-]=[\hat J^\pm,\hat N^\mp]=0 &
[\hat J^i,\hat N^i]=0\cr\ea$$
$$\Delta \hat J^3=1\otimes \hat J^3+\hat J^3\otimes 1+{\hat h\over 2}
(\hat N^3\otimes \hat J^+-\hat J^+\otimes \hat N^3)$$
$$\Delta \hat N^-=1\otimes \hat N^-+\hat N^-\otimes 1+{\hat h\over 2}
(\hat J^-\otimes \hat J^+-\hat J^+\otimes \hat J^-)$$
This is a deformation of U(iiso(2)).
\item $ ( \mu_1 , \mu_2 , \mu_3 ) = ( 0 , 1 , 0 ) $
$$\ba{ll}
[\hat J^3,\hat J^\pm ]=\pm 2\hat J^\pm & [\hat J^3,\hat N^\pm ]=\pm 2
\hat N^\pm\cr [\hat N^3,\hat N^\pm ]=[\hat N^3,\hat J^\pm ]=0 & [\hat J^+,
\hat J^-]=[\hat N^+,\hat N^-]=0\cr [\hat J^\pm,\hat N^\mp]=\pm \hat N^3 &
[\hat J^i,\hat N^i]=0\cr\ea$$
$$\Delta \hat J^3=1\otimes \hat J^3+\hat J^3\otimes 1+{\hat h\over 2}
(\hat N^3\otimes \hat J^+-\hat J^+\otimes \hat N^3)$$
This is a deformation of U(i$'$iso(2)).
\item $ ( \mu_1 , \mu_2 , \mu_3 ) = ( 0 , 0 , 0 ) $
$$\ba{ll}
[\hat J^3,\hat J^\pm ]=\pm 2\hat J^\pm & [\hat J^3,\hat N^\pm ]=\pm 2
\hat N^\pm\cr [\hat N^3,\hat N^\pm ]=[\hat N^3,\hat J^\pm ]=0 & [\hat J^+,
\hat J^-]=[\hat N^+,\hat N^-]=0\cr [\hat J^\pm,\hat N^\mp]=0 & [\hat J^i,
\hat N^i]=0\cr\ea$$
$$\Delta \hat J^3=1\otimes \hat J^3+\hat J^3\otimes 1+{\hat h\over 2}
(\hat N^3\otimes \hat J^+-\hat J^+\otimes \hat N^3)$$
This is a deformation of U(R$\oplus$(R$^4\oplus^s$so(2))).
\item $ ( \mu_1 , \mu_2 , \mu_3 ) = ( 0 , 1 , 1 ) $
$$\ba{l}
[\hat J^3,\hat J^+]=2\hat J^+\cosh({\hat h\over 2} \hat N^+)\cr
[\hat J^3,\hat J^-]=-\hat J^-\cosh({\hat h\over 2} \hat N^+)-
\cosh({\hat h\over 2} \hat N^+)\hat J^--{\hat h\over 2}\big[
\hat N^- \hat J^+\sinh({\hat h\over 2} \hat N^+)+\hat J^+
\sinh({\hat h\over 2} \hat N^+)\hat N^-\big]\cr
[\hat J^3,\hat N^+]={4\over\hat h} \sinh({\hat h\over 2} \hat N^+)\cr
[\hat J^3,\hat N^-]=-\hat N^-\cosh({\hat h\over 2} \hat N^+)-
\cosh({\hat h\over 2} \hat N^+)\hat N^-\cr
[\hat N^3,\hat N^{\pm}]=0,\qquad [\hat N^3,\hat J^{\pm}]=[\hat J^3,
\hat N^{\pm}]\cr [\hat J^+,\hat J^-]=\hat J^3,\qquad [\hat N^+,\hat N^-]=0\cr
[\hat J^{\pm},\hat N^{\mp}]=\pm \hat N^3,\qquad [\hat J^i,\hat N^i]=0\cr
  \ea$$
$$ \ba{l}
\Delta \hat J^j=e^{-\frac{\hat h}{2}\hat N^+}\otimes \hat J^j+
\hat J^j\otimes e^{\frac{\hat h}{2}\hat N^+}
-{\hat h\over 2}(e^{-\frac{\hat h}{2}\hat N^+}\hat J^+\otimes \hat N^j-
\hat N^j\otimes \hat J^+e^{\frac{\hat h}{2}\hat N^+})\cr
\Delta \hat N^j=e^{-\frac{\hat h}{2}\hat N^+}\hat J^+\otimes \hat N^j+
\hat N^j \otimes e^{\frac{\hat h}{2}\hat N^+}\qquad {\rm where\;\;}j=-,3 \cr
\gamma(X)=-e^{\hat h\hat N^+}Xe^{-\hat h\hat N^+}.\cr \ea$$
This is a deformation of U(iso(3)).
\item $ ( \mu_1 , \mu_2 , \mu_3 ) = ( 0 , 0 , 1 ) $
$$\ba{l}
[\hat J^3,\hat J^+]=2\hat J^+\cosh({\hat h\over 2} \hat N^+)\cr
[\hat J^3,\hat J^-]=-\hat J^-\cosh({\hat h\over 2} \hat N^+)-
\cosh({\hat h\over 2} \hat N^+)\hat J^--{\hat h\over 2}\big[
\hat N^- \hat J^+\sinh({\hat h\over 2} \hat N^+)+\hat J^+
\sinh({\hat h\over 2} \hat N^+)\hat N^-\big]\cr
[\hat J^3,\hat N^+]={4\over\hat h} \sinh({\hat h\over 2} \hat N^+)\cr
[\hat J^3,\hat N^-]=-\hat N^-\cosh({\hat h\over 2} \hat N^+)-
\cosh({\hat h\over 2} \hat N^+)\hat N^-\cr
[\hat N^3,\hat N^{\pm}]=0,\qquad [\hat N^3,\hat J^{\pm}]=[\hat J^3,
\hat N^{\pm}]\cr [\hat J^+,\hat J^-]=[\hat N^+,\hat N^-]=[\hat J^{\pm},
\hat N^{\mp}]=[\hat J^i,\hat N^i]=0\cr\ea$$
$$ \ba{l}
\Delta \hat J^j=e^{-\frac{\hat h}{2}\hat N^+}\otimes \hat J^j+\hat J^j
\otimes e^{\frac{\hat h}{2}\hat N^+}-{\hat h\over 2}(e^{-\frac{\hat h}{2}
\hat N^+}\hat J^+\otimes \hat N^j-\hat N^j \otimes \hat J^+
e^{\frac{\hat h}{2}\hat N^+} \cr \Delta \hat N^j=
e^{-\frac{\hat h}{2}\hat N^+}\hat J^+\otimes \hat N^j+\hat N^j \otimes
e^{\frac{\hat h}{2}\hat N^+} \qquad {\rm where\;\;}j=-,3 \cr
\gamma(X)=-e^{\hat h\hat N^+}Xe^{-\hat h\hat N^+}.\cr \ea$$
This is a deformation of U(iiso(2)).
\item $ ( \mu_1 , \mu_2 , \mu_3 ) = ( 1 , 0 , 1 ) $
$$\ba{l}
[\hat J^3,\hat J^+]={4\over\hat h} \sinh({\hat h\over 2} \hat J^+)
\cosh({\hat h\over 2} \hat N^+)\cr [\hat J^3,\hat J^-]=-\hat J^-
\cosh({\hat h\over 2} \hat J^+)\cosh({\hat h\over 2} \hat N^+)-
\cosh({\hat h\over 2} \hat J^+)\cosh({\hat h\over 2} \hat N^+)\hat J^-\cr
\hskip 1.5cm           -\hat N^-\sinh({\hat h\over 2} \hat J^+)
\sinh({\hat h\over 2} \hat N^+)-  \sinh({\hat h\over 2} \hat J^+)
\sinh({\hat h\over 2} \hat N^+)\hat N^-\cr [\hat J^3,\hat N^+]={4\over\hat h}
\sinh({\hat h\over 2} \hat N^+)\cosh({\hat h\over 2} \hat J^+)\cr
[\hat J^3,\hat N^-]=-\hat N^-\cosh({\hat h\over 2} \hat J^+)
\cosh({\hat h\over 2} \hat N^+)-\cosh({\hat h\over 2} \hat J^+)
\cosh({\hat h\over 2} \hat N^+)\hat N^-\cr\hskip 1.5cm -\hat J^-
\sinh({\hat h\over 2} \hat J^+)\sinh({\hat h\over 2} \hat N^+)-
\sinh({\hat h\over 2} \hat J^+)\sinh({\hat h\over 2} \hat N^+)\hat J^-\cr
[\hat N^3,\hat N^{\pm}]=[\hat J^3,\hat J^{\pm}],\qquad [\hat N^3,
\hat J^{\pm}]=[\hat J^3,\hat N^{\pm}]\cr [\hat J^+,\hat J^-]=[\hat N^+,
\hat N^-]=[\hat J^{\pm},\hat N^{\mp}]=[\hat J^i,\hat N^i]= 0\cr\ea$$
$$\ba{l}
\Delta \hat J^j = e^{-\frac{\hat h}{2}\hat N^+}\cosh(\frac{\hat h}{2}
\hat J^+) \otimes \hat J^j + \hat J^j \otimes \cosh(\frac{\hat h}{2}\hat J^+)
e^{\frac{\hat h}{2}\hat N^+}\cr \hskip 1cm -e^{-\frac{\hat h}{2}\hat N^+}
\sinh(\frac{\hat h}{2} \hat J^+) \otimes \hat N^j+ \hat N^j \otimes
\sinh(\frac{\hat h}{2} \hat J^+) e^{\frac{\hat h}{2}\hat N^+} \cr
\Delta \hat N^j  =  e^{-\frac{\hat h}{2}\hat N^+} \cosh(\frac{\hat h}{2}
\hat J^+)\otimes \hat N^j +\hat N^j \otimes \cosh(\frac{\hat h}{2} \hat J^+)
e^{\frac{\hat h}{2}\hat N^+} \cr \hskip 1cm -e^{-\frac{\hat h}{2}\hat N^+}
\sinh(\frac{\hat h}{2} \hat J^+) \otimes \hat J^j +\hat J^j \otimes
\sinh(\frac{\hat h}{2} \hat J^+)e^{\frac{\hat h}{2}\hat N^+}\qquad
{\rm where\;\;}j=-,3 \cr \gamma(X)=-e^{\hat h\hat N^+}Xe^{-\hat h\hat N^+}.
\cr \ea$$
This is a deformation of U(iso(2)$\oplus$iso(2)).
\end{enumerate}


\end{document}